\documentclass[number,sort&compress]{elsarticle}

\usepackage{amssymb,amsmath,amsthm,url}
\usepackage[hidelinks]{hyperref}
\usepackage{cleveref}

\crefname{equation}{}{}
\Crefname{equation}{Equation}{Equations}
\crefname{theorem}{Theorem}{Theorems}
\Crefname{theorem}{Theorem}{Theorems}
\crefname{lemma}{Lemma}{Lemmas}
\Crefname{lemma}{Lemma}{Lemmas}
\crefname{proposition}{Proposition}{Propositions}
\Crefname{proposition}{Proposition}{Propositions}
\crefname{corollary}{Corollary}{Corollaries}
\Crefname{corollary}{Corollary}{Corollaries}
\crefname{conjecture}{Conjecture}{Conjectures}
\Crefname{conjecture}{Conjecture}{Conjectures}
\crefname{section}{Section}{Sections}
\Crefname{section}{Section}{Sections}
\crefname{example}{Example}{Examples}
\Crefname{example}{Example}{Examples}

\newcommand{\Z}{\mathbb{Z}}

\newtheorem{theorem}{Theorem}[section]
\newtheorem{lemma}[theorem]{Lemma}
\newtheorem{proposition}[theorem]{Proposition}
\newtheorem{corollary}[theorem]{Corollary}

\newdefinition{example}[theorem]{Example}
\newdefinition{problem}[theorem]{Problem}

\numberwithin{equation}{section}

\newcommand{\doi}[1]{\href{http://dx.doi.org/#1}{\texttt{doi:#1}}}

\title{
On Tur\'{a}n numbers of the complete $4$-graphs
}

\author{Alexander Sidorenko}
\ead{sidorenko.ny@gmail.com}
\address{R\'{e}nyi Institute, Budapest, Hungary}

\date{\today}

\begin{document}

\begin{abstract}
The Tur\'{a}n number $T(n,\alpha+1,r)$ 
is the minimum number of edges 
in an $n$-vertex $r$-graph 
whose independence number does not exceed $\alpha$. 
For each $r\geq 2$, there exists $t_*(r)$ such that 
$T(n,\alpha+1,r) = t_*(r) \: n^r \: \alpha^{1-r} \: (1+o(1))$ 
as $\alpha / r \to\infty$ and $n / \alpha \to\infty$. 
It is known that $t_*(2) = 1/2$, 
and the conjectured value of $t_*(3)$ is $2/3$. 
We prove that $t_*(4) < 0.706335\:$.
\end{abstract}


\maketitle

\section{Introduction}

An $r$-{\it graph} is a pair $H=(V(H),E(H))$ 
where $V(H)$ is a finite set of vertices, 
and the edge set $E(H)$ is a collection of $r$-subsets of $V(H)$. 
We denote ${\mathrm v}(H)=|V(H)|$ and ${\mathrm e}(H)=|E(H)|$. 
For $A \subseteq V(H)$, 
we denote by $H \cap A$ the subgraph induced by $A$, 
that is $H \cap A = (A,E')$ 
where $E'= \{b \in E(H) : b \subseteq A\}$. 
We also denote by $H-A$ the subgraph induced by $V(H)-A$, 
that is $H-A = H \cap (V(H)-A)$. 
A subset of vertices is called {\it independent} 
if it contains no edges of $H$. 
The {\it independence number} $\alpha(H)$ 
is the largest size of an independent subset. 
The {\it transversal number} $\tau(H) = {\mathrm v}(H) - \alpha(H)$ 
is the smallest size of a subset of vertices that intersects every edge. 

The classical Tur\'{a}n number $T(n,k,r)$ 
is the minimum number of edges 
in an $n$-vertex $r$-graph $H$ with $\alpha(H) < k$. 
Consequently, $\binom{n}{r} - T(n,k,r)$ 
is the largest number of edges in an $n$-vertex $r$-graph 
that does not contain a complete subgraph on $k$ vertices.
Tables of exact values and bounds for $T(n,k,r)$ with small $n$ 
can be found in the web databases~\cite{Gordon:tables,Markstrom:tables}.

Katona, Nemetz and Simonovits \cite{Katona:1964}, 
as well as Sch\"{o}nheim \cite{Schonheim:1964}, 
showed that the ratio $T(n,k,r) / \binom{n}{r}$ is increasing with $n$, 
so the limit 
\[
  t(k,r) = \lim_{n\to\infty} T(n,k,r) / \binom{n}{r}
\]
exists and is called {\it Tur\'{a}n density}. 
It follows from the trivial recursive bound 
$T(mn,m\alpha+1,r) \leq mT(n,\alpha+1,r)$
that 
\begin{equation}\label{eq:recursive}
  t(m\alpha+1,r) \; \leq \; m^{1-r} \: t(\alpha+1,r) \: .
\end{equation}
It is convenient to use the {\it rescaled Tur\'{a}n density} 
\[
  t_*(\alpha+1,r) \; = \; \frac{\alpha^{r-1}}{r!} \: t(\alpha+1,r) \; . 
\]
While $t(k,r)$ is decreasing in parameter $k$, 
it is not clear whether the same is true for $t_*(k,r)$. 
Still, \cref{eq:recursive} yields 
$
  t_*(m\alpha+1,r) \leq t_*(\alpha+1,r)
$.
Hence, there exists the limit
\[
  t_*(r) \: = \: \lim_{k\to\infty} t_*(k,r) 
          \: = \: \inf_k t_*(k,r) \: ,
\]
and for a fixed $r$, 
when $\alpha / r \to\infty$ and $n / \alpha \to\infty$, we get 
\[ 
  T(n,\alpha+1,r) 
                  \; = \; t_*(r) \: \frac{n^r}{\alpha^{r-1}} \: (1+o(1)) 
                  \: .
\]

The exact values of $T(n,k,2)$ were found 
by Mantel \cite{Mantel:1907} for $k=3$, 
and by Tur\'{a}n \cite{Turan:1941} for all $k$. 
In particular, $t(k,2) = 1/(k-1)$ and $t_*(2)=1/2$. 
Not a single value $t(k,r)$ is known with $k>r>2$. 
It is widely believed that $t(k,3) = 4/(k-1)^2$ and $t_*(3)=2/3$. 
For arbitrary $r$, it is known that 
$t(\alpha+1,r) \leq \big((r-1) / \alpha\big)^{r-1}$ and 
$t_*(r) \leq (r-1)^{r-1} / r!$ 
(for details, see the surveys \cite{Keevash:2011,Sidorenko:1995}). 
Giraud \cite{Giraud:1990} 
discovered an elegant construction for $r=4$, $k=5$ which yields 
$t(5,4) \leq 5/16$ and $t_*(4) \leq t_*(5,4) \leq 5/6$; 
the construction was generalized 
by de Caen, Kreher and Wiseman \cite{Caen:1988}.
We will describe this construction in \cref{sec:5_16}. 

This article focuses on the case $r=4$. 
In contrast to the cases $r=2,3$, 
the values $t_*(k,4)$ vary with $k$, and in general, decrease. 
We will present new constructions 
which improve upper bounds on $t(k,4)$ for $k \geq 6$ 
and show that $t_*(4) \leq t_*(65,4) < 0.706335$\:.

A lower bound on $t_*(4)$ can be derived from the result of 
Thomass\'{e} and Yeo \cite{Thomasse:2007}. They proved that 
$\tau(H) \leq (5{\mathrm v}(H) + 4{\mathrm e}(H)) / 21$ 
holds for any $4$-graph $H$. 
This can be rewritten as 
${\mathrm e}(H) \geq 4{\mathrm v}(H) - (21/4)\alpha(H)$,
and consequently, 
$T(n,\alpha+1,4) \geq 4n - (21/4) \: \alpha$. 
In particular, $T(7m,4m+1,4) \geq 7m$.
Therefore,
$t(7m,4m+1) \geq T(7m,4m+1,4) / \binom{7m}{4} \geq 7m / \binom{7m}{4}$,
and we get a lower bound 
\[
  t_*(4) \; \geq \;
  \lim_{m\to\infty} 
    \frac{(4m)^3}{24} \: 7m \: \binom{7m}{4}^{-1}
    = \; \left(\frac{4}{7}\right)^3 
   \approx \; 0.18659 \; .
\]

Tur\'{a}n \cite{Turan:1969} conjectured that 
$t(k,r) = \left((r-1)/(k-1)\right)^{r-1}$ 
whenever $(k-1)$ is a multiple of $(r-1)$. 
While this is true for $r=2$, 
and likely to be true for $r=3$, 
our results demonstrate that it is not so for $r=4$. 
For example, we show in \cref{sec:n_7_4} that 
$t(7,4) \leq 443 / 5120 \approx 0.08652 < 1/8$.

We find 
the exact values of $T(n,\alpha+1,4)$ for $n \leq \frac{7}{4} \alpha$ 
in \cref{sec:exact}, 
and present some open problems in \cref{sec:open}.

\section{The $5/16$ construction for $T(n,5,4)$}\label{sec:5_16}

Consider two disjoint sets 
$X=\{x_1,x_2,\ldots,x_n\}$, $Y=\{y_1,y_2,\ldots,y_m\}$, 
and an $n \times m$ binary matrix $A=[a_{ij}]$. 
Let $E_{40}$ be the set of all quadruples within $X$, 
$E_{04}$ be the set of all quadruples within $Y$, 
and $E_{22}$ be the set of quadruples $\{x_i,x_j,y_k,y_l\}$ 
such that $a_{ik} + a_{il} + a_{jk} + a_{jl}$ is even.
It is easy to see that in the $4$-graph 
$H=(X \cup Y,\: E_{40} \cup E_{04} \cup E_{22})$ 
any subset of $5$ vertices contains at least one edge. 
If $n$ and $m$ are approximately equal,
and the entries of $A$ are selected randomly and independently 
with equal probability of being $0$ or $1$, 
the expected number of edges in $H$ is 
$\frac{5}{16} \binom{{\mathrm v}(H)}{4} + O({\mathrm v}(H)^3)$. 

A more specific choice of $A$ related to Hadamard matrices 
provides the best known upper bounds for $T(n,5,4)$ 
(see \cite{Sidorenko:1995}). 
The exact values of $T(n,5,4)$ were determined 
by de Caen, Kreher and Wiseman \cite{Caen:1988} 
for $n \leq 10$,
and by Markstr\"{o}m \cite{Markstrom:2009} 
for $n \leq 16$. 
The lower bound $T(17,5,4) \geq 627$ 
was also obtained in \cite{Markstrom:2009}. 
As $T(18,5,4) \geq \lceil 18 \cdot T(17,5,4) / (18-4) \rceil \geq 807$, 
we get $t(5,4) \geq 807 / \binom{18}{4} > 0.2637$ 
and $t_*(5,4) > 0.703268$.

\begin{center}
\begin{tabular}{rc}
  $n$ & $T(n,5,4)$ \\ 
\\  6 &   3
\\  7 &   7
\\  8 &  14
\\  9 &  30
\\ 10 &  50
\\ 11 &  84
\\ 12 & 126
\\ 13 & 185
\\ 14 & 259
\\ 15 & 357
\\ 16 & 476
\\ 17 & 627--644
\\ 18 & 807--828
\end{tabular}
\end{center}

\section{Expansion construction}\label{sec:expansion}

For any $4$-graph $H$, 
there is a way to augment its blow-up 
to construct a $4$-graph ${\mathcal H}$ 
with an arbitrarily large number of vertices 
and $\alpha({\mathcal H}) = \alpha(H) + 1$. 

We call a subset of vertices  $U \subseteq V$ of a $4$-graph $H=(V,E)$ 
{\it critical} if $\alpha(H-U) < \alpha(H)$. 

For each vertex $w \in V(H)$, select two critical subsets 
$I_w^0, I_w^1$ which do not contain $w$. 
These two subsets are not required to be distinct, 
but we are interested in making the independence number of 
$H - (I_w^0 \cup I_w^1)$ smaller. 
Set $d(w)=1$ if 
$\alpha\big(H - (I_w^0 \cup I_w^1)\big) \leq \alpha(H) - 2$, 
and $d(w)=0$ otherwise. 

Select a family of disjoint finite sets $\{V_w\}_{w \in V(H)}$.
First, we are going to define a $4$-graph $(V_w,E_w)$ 
in such a way that $\alpha((V_w,E_w)) = 3 + d(w)$.
 
If $d(w)=1$, 
we define $(V_w,E_w)$ as a $4$-graph described in~\cref{sec:5_16}.
Its vertex set is split in two approximately equal parts. 
Its edge set consists of all quadruples inside each of the parts 
plus approximately half of the quadruples 
which split evenly between the parts. 
We call a triple of vertices $\{a,b,c\}\subseteq V_w$ {\it even (odd)} 
if the first part contains even (odd) number of vertices from this triple. 
We set $i(a,b,c)=0$ for an even triple,
and $i(a,b,c)=1$ for an odd triple. 
Notice that if a quadruple of vertices in $(V_w,E_w)$ is not an edge, 
then it contains an even triple as well as an odd triple. 

If $d(w)=0$, we define $(V_w,E_w)$ as the complete $4$-graph. 
We set $i(a,b,c)=0$ for all triples $\{a,b,c\}\subseteq V_w$. 

An {\it expansion} of $H$ is the $4$-graph ${\mathcal H}$ 
with the vertex set $\cup_{w \in V(H)} V_w$ 
and the edge set 
\[
  E_{1111} \; \cup \;
  E_{22}   \; \cup \;
  E_{31}   \; \cup
  \bigcup_{w \in V(H)} E_w
  \: ,
\] 
where 
\[
  E_{1111} = \bigcup_{\{w,x,y,z\} \in E(H)}
    \big\{ \{a,b,c,d\} : \: 
      a \in V_w,\: b \in V_x,\: c \in V_y,\: d \in V_z \big\}
        \: ,
\]
\[
  E_{22} = \bigcup_{\{x,y\} \subseteq V(H)}
    \big\{ \{a,b,c,d\} : \: 
      \{a,b\} \subseteq V_x,\: \{c,d\} \subseteq V_y \big\}
        \: ,
\]
\[
  E_{31} = \bigcup_{w \in V(H)}
    \big\{ \{a,b,c,d\} : \: 
      \{a,b,c\} \subseteq V_w,\:  d \in V_x,\: x \in I_w^i, \: i=i(a,b,c) \big\}
        \: ,
\]

\begin{proposition}\label{th:expansion}
If ${\mathcal H}$ is an expansion of a $4$-graph $H$, 
then $\alpha({\mathcal H}) \leq \alpha(H) + 1$.
\end{proposition}

\begin{proof}[\bf{Proof}]
Let $A$ be a nonempty independent set of vertices in ${\mathcal H}$. 
We need to prove $|A| \leq  \alpha(H) + 1$. 
Let $\beta$ denote the number of sets $V_w$ which intersect $A$, 
and let $\gamma$ denote the size of the largest intersection 
$|A \cap V_w|$ over all $w \in V(H)$. 
Observe that $\beta \leq  \alpha(H)$ 
because of the way $E_{1111}$ is constructed, 
and $\gamma \leq 4$ because 
the $4$-graph $(V_w,E_w)$ does not have an independent set of size $5$. 

If there are two vertices $x,y \in V(H)$ such that $x \neq y$, 
$|A \cap V_x| \geq 2$ and $|A \cap V_y| \geq 2$, 
then $A$ contains one of the edges from $E_{22}$. 
Hence, there is at most one vertex $w \in V(H)$ such that 
$|A \cap V_w| > 1$, 
and consequently, $|A| = \beta + \gamma - 1$. 
If $\gamma \leq 2$, then $|A| \leq \alpha(H) + 2 - 1$. 
The two remaining cases are $\gamma=3,4$. 

If $\gamma = 3$ and $\beta \geq \alpha(H)$, 
then $A$ contains one of the edges from $E_{31}$. 
Hence, $\beta \leq \alpha(H)-1$ 
and $|A| \leq (\alpha(H)-1) + 3 - 1$.

In the case $\gamma = 4$, 
let $w \in V(H)$ be such that $|A \cap V_w| = 4$, 
and denote $Q = A \cap V_w$. 
As $\alpha(V_w,E_w) = 3+d(w)$, we get $d(w)=1$. 
It means that $\alpha(H - (I_w^0 \cup I_w^1)) \leq \alpha(H)-2$. 
Since quadruple $Q$ is not an edge of $(V_w,E_w)$, 
it contains both an even and and odd triple. 
If $\beta > \alpha(H)-2$, 
then $A$ contains one of the edges from $E_{31}$. 
Hence, $\beta \leq \alpha(H)-2$ 
and $|A| \leq (\alpha(H)-2) + 4 - 1$.
\end{proof}

\begin{proposition}\label{th:expansion_density}
Let $H$ be a $4$-graph with $n$ vertices. 
For each $w \in V(H)$, 
let $I_w^0, I_w^1$ be critical subsets of vertices 
that do not contain $w$, and $|I_w^0| = |I_w^1| = c(w)$. 
Let $d(w)=1$ 
if $\alpha\big(H - (I_w^0 \cup I_w^1)\big) \leq \alpha(H) - 2$, 
and $d(w)=0$ otherwise. 
Let $d=\sum_{w \in V(H)} d(w)$, 
and $c=\sum_{w \in V(H)} c(w)$.
Then 
\[
  t(\alpha(H)+2,\:4) \: \leq \:
  \big(24 {\mathrm e}(H) + 3n(n-1) + 4c + n - (11/16) d\big) \: n^{-4} \: .
\]
\end{proposition}

\begin{proof}[\bf{Proof}]
Let ${\mathcal H_N}$ be an expansion of $H$ where 
the sets $V_w$ are of size $N$. 
Then 
\[
  |E_{1111}| = {\mathrm e}(H) N^4 , \;\;\;\;
  |E_{22}|   = \binom{n}{2} \binom{N}{2}^2 , \;\;\;\;
  |E_{31}|   =  \sum_{w \in V(H)} c(w) N \binom{N}{3},
\]
\[
  \big| \bigcup_{w \in V(H)} E_w \big| \: = \:
  \left((n-d) \cdot 1 + d \cdot \frac{5}{16}\right) \binom{N}{4} + O(N^3)
                   \;\;\mbox{as}\;\; N\to\infty \; .
\]
Hence,
\[
  {\mathrm e}({\mathcal H_N}) =
  \big[24 {\mathrm e}(H) + 3n(n-1) + 4c + n - (11/16)d \big] \binom{N}{4}
  + O(N^3) \, .
\]
By \cref{th:expansion}, $\:\alpha({\mathcal H_N}) \leq \alpha(H)+1$. 
Therefore,
\begin{eqnarray*}
  t(\alpha(H)+2,4) & \leq & 
  \lim_{N\to\infty} 
    {\mathrm e}({\mathcal H_N}) \left/ \binom{nN}{4} \right. \\
  & = & 
  \big(24 {\mathrm e}(H) + 3n(n-1) + 4c + n - (11/16)d \big) \: n^{-4} \: .
\end{eqnarray*}
\end{proof}

\begin{example}\label{ex:1}
Let $H$ be a $4$-graph with $8$ vertices $1,2,\ldots,8$ 
and $14$ edges: 
$\{1,2,3,4\}$, 
$\{5,6,7,8\}$, 
$\{1,2,5,6\}$, 
$\{3,4,7,8\}$, 
$\{1,2,7,8\}$, 
$\{3,4,5,6\}$, 
$\{1,3,5,7\}$, 
$\{2,4,6,8\}$, 
$\{1,3,6,8\}$, 
$\{2,4,5,7\}$, 
$\{1,4,5,8\}$, 
$\{2,3,6,7\}$, 
$\{1,4,6,7\}$, 
$\{2,3,5,8\}$. 
Then $\alpha(H)=4$, and each edge is a critical subset. 
Let $I_i^0,I_i^1$ be two different edges of $H$ which do not contain $i$. 
Then $\alpha(H-(I_i^0 \cup I_i^1)) = 2$. 
By \cref{th:expansion_density}, we get $t(6,4) \leq 1269/8192$ 
and $t_*(6,4) \leq 5^3 \cdot 423 \cdot 2^{-16} < 0.80681\:$.
\end{example}

In the proof of \cref{th:expansion_density}, 
we used the sets $V_w$ of equal size. 
Often, it is advantageous to select 
critical subsets $I_i^0,I_i^1$ in an asymmetric way 
and exploit this asymmetry 
by optimizing the sizes of sets $V_w$.

\begin{example}\label{ex:2}
We will use $4$-graph $H$ from Example~\ref{ex:1} 
and select the critical subsets as follows: 
$I_1^0 = I_2^0 = I_3^0 = I_4^0 = \{5,6,7,8\}$, 
$I_1^1 = I_2^1 = I_5^0 = I_6^0 = \{3,4,7,8\}$, 
$I_3^1 = I_4^1 = I_5^1 = I_6^1 = \{1,2,7,8\}$, 
$I_7^0 = \{2,4,6,8\}$, 
$I_7^1 = \{1,3,6,8\}$, 
$I_8^0 = \{2,3,6,7\}$, 
$I_8^1 = \{1,4,6,7\}$. 
Let $N$ be the number of vertices in the expansion. 
If $|V_i| = x_i N + O(1)$ where
$x_1=x_2=x_3=x_4=0.13387$, 
$x_5=0.13639$, 
$x_6=0.13085$, 
$x_7=x_8=0.09684$, 
then the number of edges in the expansion is 
$(a + O(1/N)) 5^{-3} N^4$ 
as $N\rightarrow\infty$, where $a < 0.80261\:$. 
Therefore, $t_*(6,4) < 0.80261\:$.
\end{example}

\section{Circular construction}\label{sec:circular}

To define a circular construction, we start with 
a circularly ordered sequence of $4$-graphs $\{G_i\}_{i\in\Z_m}$ 
where $m \geq 2$. 
The vertex set of each $G_i$ is partitioned into two disjoint subsets: 
$V(G_i) = V_i' \cup V_i''$. 
We denote $W_i = V_i' \times V_{i+1}''$, 
and for $x \in V_i'$, denote $W_{i,x} = \{(x,y): y \in V_i''\}$.
For a subset $A \subseteq (W_0 \cup W_1 \cup \ldots \cup W_{m-1})$,
let $\beta_i(A)$ denote the number of sets $W_{i,x}$ that intersect $A$. 

We define functions 
$f_i:\: W_i \to V_i'$, 
$g_i:\: W_i \to V_{i+1}''$, 
and $h_i: (W_i \cup W_{i+1}) \to V_{i+1}''$ 
as follows. 
For $w=(x,y) \in W_i$, where $x \in V_i'$ and $y \in V_{i+1}''$, 
we set $f_i(w)=x$ and $g_i(w)=y$. 
For $w \in (W_i \cup W_{i+1})$, 
we set $h_i(w)=g_i(w)$ if $w \in W_i$, 
and $h_i(w)=f_{i+1}(w)$ if $w \in W_{i+1}$. 

Let $E_i^1$ denote the set of quadruples 
$\{w_1,w_2,w_3,w_4\} \subseteq W_i \cup W_{i+1}$ such that 
$\beta_i(\{w_1,w_2,w_3,w_4\})=1$, 
and $h_i(w_1),h_i(w_2),h_i(w_3),h_i(w_4)$ 
form an edge in $G_{i+1}$. 

Let $E_i^2$ denote the set of quadruples 
$\{w_1,w_2,w_3,w_4\} \subseteq W_i$ such that $w_j = (x_j,y_j)$, 
$\;x_1,x_2,x_3,x_4 \in V_i'$, $\;x_1=x_2$, $\:x_3=x_4 \neq x_1$. 

Let $E_i^4$ denote the set of quadruples 
$\{w_1,w_2,w_3,w_4\} \subseteq W_i$ such that $w_j = (x_j,y_j)$, 
where 
$\;x_1,x_2,x_3,x_4 \in V_i'$ are pairwise distinct 
and form an edge in $G_i$. 

The {\it circular construction} $C_m[G_0,G_1,\ldots,G_{m-1}]$ 
is a $4$-graph with the vertex 
set $W_0 \cup W_1 \cup \ldots \cup W_{m-1}$ 
and the edge set
\[
  \bigcup_{i\in\Z_m} \big( E_i^1 \cup E_i^2 \cup E_i^4 \big)
  \: .
\]
As 
\begin{eqnarray*}
  |E_i^1| & = & |V_i'| \: \sum_{j=0}^3 
    \left|\{
      e \in E(G_{i+1}):\: |e \cap V_{i+1}'|=j
    \}\right|
    \cdot |V_{i+2}''|^j
    ,
\\
  |E_i^2| & = & \binom{|V_i'|}{2} \binom{|V_{i+1}''|}{2}^2
    ,
\\
  |E_i^4| & = & {\rm e}(G_i \cap V_i') \cdot |V_{i+1}''|^4
    ,
\end{eqnarray*}
we get 
\begin{eqnarray*}
  {\rm v}\left(C_m[G_0,G_1,\ldots,G_{m-1}]\right) & = &
  \sum_{i\in\Z_m} |V_i'| \cdot |V_{i+1}''| \: ,
    \;\;\;\;\;\;\;\;\;\;\;\;\;\;\;\;\;\;\;\;\;\;\;\;\;
  \\ \\
  {\rm e}\left(C_m[G_0,G_1,\ldots,G_{m-1}]\right) & = &
\end{eqnarray*}
\begin{eqnarray}\label{eq:circ}
  & = &
  \sum_{i\in\Z_m} 
    |V_{i-1}'| \: \sum_{j=0}^3 \left|\{e \in E(G_i):\:
      |e \cap V_i'|=j\}\right| \cdot |V_{i+1}''|^j
    \nonumber \\
  & + & \sum_{i\in\Z_m}
    \binom{|V_i'|}{2} \binom{|V_{i+1}''|}{2}^2
  \; + \; \sum_{i\in\Z_m}
    {\rm e}(G_i \cap V_i') \cdot |V_{i+1}''|^4
    .
\end{eqnarray}

\begin{theorem}\label{th:circ1}
If $\alpha(G_i) = \alpha_i + 1$ 
and $\alpha(G \cap V_i') \leq \alpha_i$ 
for every $i \in \Z_m$, 
then $\alpha(C_m[G_0,G_1,\ldots,G_{m-1}]) \leq \sum_{i\in\Z_m} \alpha_i$.
\end{theorem}

\begin{proof}[\bf{Proof}]
Let $A$ be an independent subset of vertices in 
$C_m[G_0,G_1,\ldots,G_{m-1}]$. 
For $i\in\Z_m$, let $\beta_i=\beta_i(A)$, 
so $\beta_i$ is the number of sets $V_{ix}$ that intersect $A$. 
Let $\gamma_i$ denote the size of the largest of these intersections. 
Set $\chi_i=0$ if $A \cap W_i = \emptyset$, 
and $\chi_i=1$ otherwise.

If there are two different $x,x' \in V_i'$ such that 
$|A \cap V_{ix}| \geq 2$ and $|A \cap V_{ix'}| \geq 2$, 
then $A$ contains a quadruple from $E_i^2$. 
Hence, if $A \cap W_i \neq \emptyset$, 
then $|A \cap W_i| = \beta_i + \gamma_i - 1$. 
In any case, $|A \cap W_i| = \beta_i + \gamma_i - \chi_i$. 

If $A \cap W_i \neq \emptyset$ 
and $\beta_{i+1} + \gamma_i > \alpha_{i+1} + 1$, 
then $A$ contains a quadruple from $E_i^1$.
If $A \cap W_i = \emptyset$, then $\gamma_i = 0$. 
If $\beta_{i+1} > \alpha_{i+1}$, 
then $A$ contains a quadruple from $E_{i+1}^4$. 
Hence, $\beta_{i+1} + \gamma_i \leq \alpha_{i+1} + \chi_i$.
So we get
\[
  |A| = \sum_{i\in\Z_m}  |A \cap W_i| 
      \: =  \sum_{i\in\Z_m} (\beta_i + \gamma_i - \chi_i) 
      \: = \sum_{i\in\Z_m} (\beta_{i+1} + \gamma_i - \chi_i) 
      \: \leq \: \sum_{i\in\Z_m} \alpha_{i+1} .
\]
\end{proof}

\begin{theorem}\label{th:circ2}
Let $A$ be an independent set 
in $C_m[G_0,G_1,\ldots,G_{m-1}]$. 
For $i\in\Z_m$, set $\chi_i=0$ if $A \cap W_i = \emptyset$, 
and $\chi_i=1$ otherwise. 
If, in addition to the conditions of \cref{th:circ1},
$\alpha(G \cap V_i'') \leq \alpha_i$, then 
$|A| \leq \sum_{i\in\Z_m} (\alpha_i - 1 + \chi_i)$.
\end{theorem}

\begin{proof}[\bf{Proof}]
The inequality $\beta_{i+1} + \gamma_i \leq \alpha_{i+1} + \chi_i$ 
in the proof of \cref{th:circ1} 
can be strengthen to 
$\beta_{i+1} + \gamma_i \leq \alpha_{i+1} + \chi_i + \chi_{i+1} - 1$. 
Indeed, if $\chi_{i+1}=0$, then $\beta_{i+1}=0$. 
If $\gamma_i > \alpha_{i+1} + \chi_i - 1$, then $\chi_i = 1$ 
and $\gamma_i > \alpha_{i+1}$. 
In this case, $A$ would contain a quadruple from $E_i^1$. 
Hence, 
$\beta_{i+1} + \gamma_i \leq \alpha_{i+1} + \chi_i + \chi_{i+1} - 1$,
and we get
\begin{eqnarray*}
  |A| & = & \sum_{i\in\Z_m}  |A \cap W_i| 
      \: =  \sum_{i\in\Z_m} (\beta_i + \gamma_i - \chi_i)
      \\ & = &
      \sum_{i\in\Z_m} (\beta_{i+1} + \gamma_i - \chi_{i+1}) 
      \: \leq \sum_{i\in\Z_m} (\alpha_i - 1 + \chi_i) .
\end{eqnarray*}
\end{proof}

\begin{example}\label{ex:4x4}
Let $G$ be a $4$-graph with the vertex set $\Z_2^3$ 
where pairwise distinct vertices $z_1,z_2,z_3,z_4$ 
form an edge if $z_1+z_2+z_3+z_4=0$. 
We partition its vertex set $V(G) = V_0 \cup V_1$ where 
$V_i$ contains vectors with first entry equal to $i\in\Z_2$. 
Then $\alpha(G)=4$ and $\alpha(G \cap V_i)=3$. 
Set $G_m = C_m[G,G,\ldots,G]$. 
In this case, 
$|W_i| = 4 \cdot 4 = 16$, 
$|E_i^1| = 772$, 
$|E_i^2| = 216$, 
$|E_i^4| = 256$, 
so we get ${\rm v}(G_m) = 16m$, ${\rm e}(G_m) = 1244m$, 
and by \cref{th:circ1}, $\alpha(G_m) \leq 3m$. 
(We skip the calculations that lead to ${\mathrm e}(H_m)=1244m$, 
because we are going to revisit this example in \cref{sec:t_65_4}.) 
Hence, $T(16m,3m+1,4) \leq 1244m$ for $m \geq 2$. 
By \cref{th:circ2}, 
we get $\alpha(G_m - W_i) \leq 3m-1$, 
and $\alpha(G_m - (W_i \cup W_j)) \leq 3m-2$ when $i \neq j$. 
We are going to build an extension of $G_m$ 
to get an upper bound on $t_*(3m+2,4)$. 
When $m \geq 3$, for any vertex $w \in W_k$, we can select 
$I_w^0 = W_i$, $I_w^1 = W_j$, where $i,j \neq k$, $i \neq j$. 
Then we get $d(w)=1$ for each $w$. 
We use \cref{th:expansion_density} 
with $n=16m$, ${\mathrm e}(H_m)=1244m$, $\alpha(H_m) = 3m$, 
$c(w)=16$ and $d(w)=1$ for all $w \in V(H_m)$ 
to derive
\begin{eqnarray*}
    t(3m+2,4) & \leq & (768m + 30,837) / (65,536 \: m^3)
      \: , \\
  t_*(3m+2,4) & \leq & (768m + 30,837) (3m+1)^3 / (1,572,864 \: m^3) \: .
\end{eqnarray*}
When $m=7$, we get $t_*(23,4) < 0.714739 \:$. 
\end{example}

\section{$T(n,6,4)$}\label{sec:n_6_4}

The best asymptotical bound we know is 
$t_*(6,4) < 0.80261$ from Example~\ref{ex:2}. 
We describe below a few elegant constructions for small $n$. 

\subsection{$n=9,10$}

Notice that any simple graph with $5$ vertices and $6$ edges 
contains either a vertex of degree $4$, or a $4$-cycle. 
Consider a $4$-graph $H$ 
whose vertices are the edges of complete $5$-vertex graph $K_5$, 
and the edges of $H$ correspond to 
$4$-arm stars and $4$-cycles in $K_5$. 
The number of edges in $H$ is $5+15=20$, so $T(10,6,4) \leq 20$. 
By removing a vertex, 
we get a $4$-graph with $9$ vertices and $12$ edges.

\subsection{$n=11,12$}

Color the edges of $K_6$ with 5 colors, 
so that edges of each color form a perfect matching. 
Take two copies of this $K_6$. 
Consider a $4$-graph whose vertex-set  
consists of the vertices of the two $K_6$, 
and each edge is a union of an edge from one $K_6$ 
and an edge from the other $K_6$ that have the same color. 
So far, we have selected $5 \cdot 3^2 = 45$ edges. 
In each $K_6$, select $3$ quadruples of vertices 
in such a way that 
every vertex belongs to exactly two of these quadruples. 
Add these $6$ quadruples as edges to our $4$-graph.
It is easy to check that now any $6$ vertices contain an edge, so 
$T(12,6,4) \leq 51$. 
This construction is not unique, 
as we can choose the $6$ quadruples in $7$ nonequivalent ways,
depending of the colors of their pairwise intersections. 
According to \cite{Markstrom:2009}, 
there are exactly $7$ nonisomorphic $4$-graphs 
with $12$ vertices, $51$ edges and independence number $5$. 
Hence, our construction captures all of them. 
By removing a vertex in any of them, 
we get a $4$-graph with $11$ vertices and $34$ edges.

\subsection{$n=14,15,16$}\label{sec:16_6_4}

Consider a $4$-graph $G$ with the vertex set $A \cup B$ 
where $A$ and $B$ are disjoint copies of $\Z_2^2 \oplus \Z_2$. 
We will denote elements of $A$ and $B$ by $(x,a)$ 
where $x\in\Z_2^2$ and $a\in\Z_2$. 
The edges of $G$ are quadruples of distinct vertices 
$v_i=(x_i,a_i)$ ($i=1,2,3,4$) where either 

({\it i}) $v_1,v_2,v_3,v_4 \in A$, \;$x_1+x_2+x_3+x_4=0$; or 

({\it ii}) $v_1,v_2,v_3,v_4 \in B$, \;$a_1+a_2+a_3+a_4=0$; or 

({\it iii}) $v_1,v_2 \in A$, $v_3,v_4 \in B$, 
            $x_1=x_2$, $a_3 \neq a_4$; or 

({\it iv}) $v_1,v_2 \in A$, $v_3,v_4 \in B$, 
           \;$a_3=a_4$, \;$x_1+x_2+x_3+x_4=0$. 

\noindent
One may check that any $6$ vertices contain an edge, so 
$T(16,6,4) \leq {\rm e}(G) = 220$. 
By removing a vertex from set $B$, 
we get a $4$-graph with $15$ vertices and $161$ edges.
By removing vertices $(x,0)$ and $(x,1)$ from set $B$, 
we get a $4$-graph with $14$ vertices and $115$ edges.

\begin{center}
\begin{tabular}{rcl}
  $n$ & $T(n,6,4)$ & reference \\ 
\\  7 &  3 & \cite{Markstrom:tables}
\\  8 &  6 & \cite{Markstrom:tables}
\\  9 & 12 & \cite{Markstrom:tables}
\\ 10 & 20 & \cite{Markstrom:tables}
\\ 11 & 34 & \cite{Markstrom:tables}
\\ 12 & 51 & \cite{Markstrom:tables}
\\ 13 &  74--79 & \cite{Gordon:tables}
\\ 14 & 104--115 & \cref{sec:16_6_4}
\\ 15 & 142--161 & \cref{sec:16_6_4}
\\ 16 & 190--220 & \cref{sec:16_6_4}
\end{tabular}
\end{center}

\section{$T(n,7,4)$}\label{sec:n_7_4}

\begin{proposition}\label{th:rainbow}
\begin{eqnarray*}
  t(7,4) & \leq &  443 / 5120 \; ,
\\
t_*(7,4) & \leq & 3987 / 5120 \;
  < \; 0.778711 \; .
\end{eqnarray*}
\end{proposition}

\begin{proof}[\bf{Proof}]
Consider a space $V^k$ 
of $k$-dimensional vectors over $\Z_2^2$. 
For each pair of distinct vectors 
$x=(x_1,x_2,\ldots,x_k)$, 
$y=(y_1,y_2,\ldots,y_k)$,
we select a certain nonzero element of $\Z_2^2$ 
to serve as its color. 
Let $i=i(x,y)$ be the smallest index such that $x_i \neq y_i$. 
We define the color of $\{x,y\}$ as $c(x,y) = x_i + y_i$. 
It is known (see \cite{Balogh:2017}) 
that this coloring produces 
the maximum possible number of rainbow triples, 
which is 
$\left(\frac{2}{5}+O(4^{-k})\right)\binom{4^k}{3}$. 

Observe that $\{x,y,z\}$ is a rainbow triple 
(that is, $c(x,y)$, $c(x,z)$, $c(y,z)$ are distinct) 
if and only if $i(x,y)=i(x,z)=i(y,z)$. 
We say that $x$ is the {\it apex} of $\{x,y,z\}$ 
if $i(x,y)=i(x,z)<i(y,z)$. 
Then every non-rainbow triple has exactly one apex. 

Next we are going to construct a $4$-graph 
whose vertex set is the union
of disjoint copies of $V'=V^1$ and $V''=V^k$, 
and the edge-set is the union of four families of quadruples 
on $V' \cup V''$ denoted below as $E_0,E_1,E_2,E_4$. 

Let $E_4$ consist of a single quadruple 
that is the set of elements of $V'$. 

Let $E_2$ be the family of such quadruples $Q$ 
that $|Q \cap V'| = |Q \cap V''| = 2$ 
and the color of $Q \cap V'$ in $V'$ 
is the same as the color of $Q \cap V''$ in $V''$. 
For any pair in $V''$, 
there are exactly two choices of a pair in $V'$ with the same color. 
Thus, $|E_2| = 2\binom{|V''|}{2} = 4^k (4^k - 1)$. 

We partition $V'$ into two parts, 
$X'=\{(0,0),(0,1)\}$ and $Y'=\{(1,0),(1,1)\}$. 
For a non-rainbow triple $T=\{x,y,z\} \subset V''$ with apex $x$, 
we define a pair $p(T)$ in $V'$ as follows. 
If $x_1,y_1 \in X'$ or $x_1,y_1 \in Y'$, 
we set $p(T)=\{y_1,y_1+c(y,z)\}$, 
and otherwise set $p(T)=V' \backslash \{y_1,y_1+c(y,z)\}$. 
(Notice that $y_1=z_1$ for non-apex elements $y$ and $z$.) 
Let $E_1$ be the family of quadruples $\{x,y,z,w\}$ such that 
$T=\{x,y,z\} \subset V''$ is a non-rainbow triple 
and $w \in p(T)$.
Then the size of $E_1$ is 
twice the number of non-rainbow triples, that is 
$|E_1| = \left(\frac{6}{5}+O(4^{-k})\right)\binom{4^k}{3}$. 

Similarly to the partition $(X',Y')$ of $V'$, 
we partition $V''$ into $X''$ and $Y''$. 
Namely, $x \in V''$ belongs to $X''$ if and only if $x_1 \in X'$. 
Let $E_0$ be the family of quadruples $Q \subseteq V''$ such that 
$|Q \cap X''|$ is even. 
Then $|E_0| = 2\binom{\frac{1}{2}4^k}{4} 
+ \binom{\frac{1}{2}4^k}{2}^{\!_2}$. 

Let $H_k$ be the $4$-graph with vertex set $V' \cup V''$ 
and edge set $E_0 \cup E_1 \cup E_2 \cup E_4$. 
We claim that $\alpha(H_k)=4$. 
Indeed, suppose $A \subset V' \cup V''$, $\:|A|=5$. 
To show that $A$ contains an edge of $H_k$, 
we consider $5$ separate cases, depending on the size of $A \cap V'$. 

If $|A \cap V'|=4$, 
then $A$ contains $V'$, and $V' \in E_4$. 

If $|A \cap V'|=3$ and $|A \cap V''|=2$, 
then $A \cap V'$ is a rainbow triple in $V'$, 
thus it contains a pair of the same color as $A \cap V''$ in $V''$, 
so $A$ contains a quadruple from $E_2$. 

If $|A \cap V'|=2$ and $|A \cap V''|=3$, 
denote $T = A \cap V''$. 
If $T$ contains a pair of the same color as $A \cap V'$, 
then $A$ contains a quadruple from $E_2$. 
If $T$ does not contain a pair of that color, 
then $T$ is not a rainbow triple. 
Let $T=\{x,y,z\}$ where $x$ is the apex. 
As the color of $\{y,z\}$ differs from the color of $A \cap V'$, 
then $(A \cap V') \cap p(T) \neq \emptyset$.
Let $w \in (A \cap V') \cap p(T)$. 
Then $T \cup \{w\}$ belongs to $E_1$. 

In the case $|A \cap V'|=1$ and $|A \cap V''|=4$, 
if $|A \cap X''|$ is even, then $A \cap V''$ belongs to $E_0$. 
Hence, we may assume that $|A \cap X''|$ is odd. 
Furthermore, without loss of generality, we may assume that 
$|A \cap X''|=3$ and $|A \cap Y''|=1$. 
Let $T = A \cap X'' = \{x,y,z\}$, $\:A \cap Y'' = \{u\}$,
$\:A \cap V' =\{w\}$. 
If $T$ is a rainbow triple, then $x_1=y_1=z_1$, 
and $u$ is the apex in triples 
$\{u,x,y\}$, $\:\{u,x,z\}$, $\:\{u,y,z\}$. 
Then the union of 
$  p(\{u,x,y\}) = \{x_1,x_1+c(x,y)\}$, 
$\:p(\{u,x,z\}) = \{x_1,x_1+c(x,z)\}$, and 
$  p(\{u,y,z\}) = \{x_1,x_1+c(y,z)\}$ 
covers $V'$ entirely, 
so one of the quadruples 
$\{u,x,y,w\}$, $\:\{u,x,z,w\}$, $\:\{u,y,z,w\}$ 
belongs to $E_1$. 
Thus, we may assume that $T$ is a non-rainbow triple. 
For definiteness, let $x$ be its apex. 
Since $p(\{x,y,z\}) = \{y_1,y_1+c(y,z)\}$ 
and $p(\{u,y,z\}) = V' \backslash \{y_1,y_1+c(y,z)\}$, 
then $p(\{x,y,z\}) \cup p(\{u,y,z\}) = V'$. 
Hence, one of the quadruples 
$\{x,y,z,w\}$ and $\{u,y,z,w\}$ 
belongs to $E_1$. 

If $|A \cap V''|=5$, then $A \subset V''$, 
so one of the quadruples in $A$ 
has an even size intersection with $X''$, 
and hence, belongs to $E_0$. 

Now we use the circular construction $H_k^m=C[H_k,H_k,\ldots,H_k]$ 
with $m$ disjoint copies of $H_k$. 
Then ${\rm v}(H_k^m) = m \: 4^{k+1}$. 
As $\alpha(H_k)=4$, we get by \cref{th:circ1}, $\alpha(H_k^m)=3m+1$. 
As ${\rm e}(H_k \cap V')=1$ 
and $|\{ e \in E(H_k): |e \cap V'|=j \}|$ is equal to
$|E_j|$ for $j=0,1,2$ and equal to $0$ for $j=3$, 
we get by~(\ref{eq:circ}),
\begin{eqnarray*}
  {\rm e}(H_k^m) & = & m \left[
    4 \left( |E_0| + |E_1| \cdot 4^k + |E_2| \cdot 4^{2k}
    \right) + 6 \binom{4^k}{2}^2 + 4^{4k}
  \right] 
  \\ & = & m \left[
    4 \left( \frac{1}{2} + \frac{24}{5} + 24
    \right) + 6 \cdot 6 + 24 + O \left( 4^{-k} \right)
  \right] \binom{4^k}{4}
  \\ & = & m \: \binom{4^k}{4} \left[ 
    \frac{886}{5}
    + O \left( 4^{-k} \right)
  \right] 
  \\ & = & m \: \binom{4^{k+1}}{4} \left[ 
    \frac{443}{640} 
    + O \left( 4^{-k} \right)
  \right] 
    \; .
\end{eqnarray*}
Then 
\[
  t(3m+1,4) \; \leq \;
  \lim_{k\to\infty} {\rm e}(H_k^m) 
    \left/ \binom{{\rm v}(H_k^m)}{4} \right.
  \, = \: \frac{443}{640} \, m^{-3} .
\]
When $m=2$, we get 
$t(7,4) \leq \frac{443}{5120}$ 
and 
$t_*(7,4) = \frac{6^3}{4!} t(7,4) = 9 \cdot t(7,4)$.
\end{proof}

\vspace{2mm}
\begin{center}
\begin{tabular}{rcl}
  $n$ & $T(n,7,4)$ & reference \\ 
\\  8 &  2 & \cite{Markstrom:tables}
\\  9 &  5 & \cite{Markstrom:tables}
\\ 10 & 10 & \cite{Markstrom:tables}
\\ 11 & 17 & \cite{Markstrom:tables}
\\ 12 & 26 & \cite{Markstrom:tables}
\\ 13 & 38--39 & \cite{Gordon:tables}
\\ 14 & 54--56 & \cite{Gordon:tables}
\\ 15 & 74--80 & \cite{Gordon:tables}
\\ 16 & 99--108 & $H_2$ from Example~\ref{ex:4x4}
\end{tabular}
\end{center}

\section{$\Z_m \oplus \Z_2^6$ construction}\label{sec:t_65_4}

In this section, we will prove that $t_*(4) < 0.706335\:$.

\medskip

We denote elements of $\Z_m \oplus \Z_2^6$ 
by $(i,x,y,z)$ with $i\in\Z_m$ and $x,y,z\in\Z_2^2$.
For $m\geq 4$, let $H_m$ be a $4$-graph with vertex set 
$\Z_m \oplus \Z_2^6$ 
where distinct vertices $v_1,v_2,v_3,v_4$ form an edge if either 
\begin{enumerate}
\item\label{it:101}
$v_t=(i,x,y,z_t)$ for $t=1,2,3,4$; or
\item\label{it:102}
$v_t=(i,x_t,y_t,z_t)$ for $t=1,2,3,4$,
where $x_1+x_2+x_3+x_4=0$, 
and $(x_k,y_k) \neq (x_l,y_l)$ for $k \neq l$; or
\item\label{it:103}
$v_1=(i,x_1,y_1',z_1')$, 
$v_2=(i,x_1,y_1'',z_1'')$, 
$v_3=(i+1,x_2',y_2',z_2')$, 
$v_4=(i+1,\linebreak x_2'',y_2'',z_2'')$, 
where $y_1' \neq y_1''$ and $y_1'+y_1''+x_2'+x_2''=0$; or
\item
$v_1=(i,x_1,y_1,z_1')$, 
$v_2=(i,x_1,y_1,z_1'')$, 
where $z_1' \neq z_1''$ and either
\begin{enumerate}
\item\label{it:104a}
$v_3=(i,x_2,y_2,z_2')$, 
$v_4=(i,x_2,y_2,z_2'')$, 
where $(x_1,y_1) \neq (x_2,y_2)$ and $z_2' \neq z_2''$; or
\item\label{it:104b}
$v_3=(i+2,x_2,y_2',z_2')$, 
$v_4=(i+2,x_2,y_2'',z_2'')$, 
where $y_2' \neq y_2''$; or
\item\label{it:104c}
$v_3=(i+3,x_2',y_2',z_2')$, 
$v_4\!=\!(i+3,x_2'',y_2'',z_2'')$, 
where $z_1'+z_1''+x_2'+x_2''\!=\!0$. 
\end{enumerate}
\end{enumerate}
Let $A$ be a subset of the vertex set of $H_m$. 
For $i \in \Z_m$ and  $x,y \in \Z_2^2$, we denote
\[
  A_1(i) \: = \: 
  \{ x \in \Z_2^2 \: | \; \exists y,z \in \Z_2^2 : \: (i,x,y,z) \in A \},
\]
\[
  A_2(i,x) \: = \: 
  \{ y \in \Z_2^2 \: | \; \exists z \in \Z_2^2 : \: (i,x,y,z) \in A \},
\]
\[
  A_3(i,x,y) \: = \: 
  \{ z \in \Z_2^2 \: | \: (i,x,y,z) \in A \},
\]
\[
  \varepsilon_i(A) \: =
  \sum_{x,y\in\Z_2^2} \max \{0,\: |A_3(i,x,y)|-1\} .
\]
Let $\chi[B]$ be the indicator of condition $B$ 
(that is, $\chi[B]$ is $1$ if $B$ is true, and $0$ otherwise).

\begin{lemma}\label{th:5dim_L1}
If $A$ is an independent set of vertices in $H_m$, then
\begin{multline*}
  \varepsilon_i(A) \: + \: 
  \max_x |A_2(i+2,x)| \: + \:
  |A_1(i+3)| \\
  \leq \: 2 \: + \: 
  \chi[A_1(i+2) \neq \emptyset] \: + \: 
  \chi[A_1(i+3) \neq \emptyset] .
\end{multline*}
\end{lemma}

\begin{proof}[\bf{Proof}]
{\it Case 1:} $\varepsilon_i(A) \! > \! 0$. 
It means that there are $x'\!,y'$ such that $|A_3(i,x',y')| \!\geq\! 2$. 
Hence, there are $z' \neq z''$ such that 
$v'=(i,x',y',z')$ and $v''=(i,x',y',z'')$ belong to $A$. 
If there is $x$ such that $|A_2(i+2,x)| \geq 2$, 
then there are $y_1 \neq y_2$ and $z_1,z_2$ such that 
$v_1=(i+2,x,y_1,z_1)$ and $v_2=(i+2,x,y_2,z_2)$ belong to $A$. 
Then $\{v',v'',v_1,v_2\}$ is an edge of type \ref{it:104b}. 
Thus, we may assume that $|A_2(i+2,x)| \leq 1$ for all $x$. 
If $A_1(i+2) = \emptyset$ then $\max_x |A_2(i+2,x)| = 0$. 
Hence, $\max_x |A_2(i+2,x)| \leq \chi[A_1(i+2) \neq \emptyset]$. 
It remains to prove that 
$\varepsilon_i(A) + |A_1(i+3)| \leq 2 + \chi[A_1(i+3) \neq \emptyset]$.

If $|A_1(i+3)| \geq 3$, 
then there are $v_j=(i+3,x_j,y_j,z_j) \in A$ $(\:j=1,2,3)$ 
where $x_1,x_2,x_3$ are pairwise distinct. 
There exist $j,k\in\{1,2,3\}$ such that $j \neq k$ and $x_j+x_k=z'+z''$. 
Then $\{v',v'',v_j,v_k\}$ is an edge of type \ref{it:104c}. 
Thus, we may assume that $|A_1(i+3)| \leq 2$. 

If there is $(x'',y'') \neq (x',y')$ 
such that $|A_3(i,x'',y'')| \geq 2$, 
then there exist $z''' \neq z''''$ such that 
$v'''=(i,x'',y'',z'''')$ and $v''''=(i,x'',y'',z'''')$ belong to $A$. 
In this case, $\{v',v'',v''',v'''\}$ is an edge of type \ref{it:104a}. 
Thus, we may assume that $|A_3(i,x'',y'')| \leq 1$ 
for any $(x'',y'') \neq (x',y')$. 
It means that $\varepsilon_i(A) = |A_3(i,x',y')| - 1$. 

If $\varepsilon_i(A) \geq 3$, 
then $|A_3(i,x',y')|=4$,  
so there are $v_j=(i,x',y',z_j) \in A$ $\:(j=1,2,3,4)$ 
where $z_1,z_2,z_3,z_4$ are pairwise distinct. 
In this case, $\{v_1,v_2,v_3,v_4\}$ is an edge of type \ref{it:101}. 
Thus, we may assume that  $\varepsilon_i(A) \leq 2$. 

If $|A_1(i+3)|=2$,
then there are $v_j=(i+3,x_j,y_j,z_j) \in A$ $\:(j=1,2)$ 
where $x_1 \neq x_2$. 
If additionally to that, $\varepsilon_i(A) \geq 2$, 
then $|A_3(i,x',y')| \geq 3$, which means that  
there exists $z''' \neq z',z''$ such that $v'''=(i,x',y',z''') \in A$. 
We may find among $v',v'',v'''$ two vertices 
which together with $v_1,v_2$ form an edge of type \ref{it:104c}. 
Thus, if $|A_1(i+3)|=2$ then  $\varepsilon_i(A) \leq 1$. 
Hence, $\varepsilon_i(A) + |A_1(i+3)| \leq 3$ 
where the equality is possible only if $|A_1(i+3)| > 0$. 
Therefore, 
$\varepsilon_i(A) + |A_1(i+3)| \leq 2 + \chi[A_1(i+3) \neq \emptyset]$.

{\it Case 2:} $\varepsilon_i(A)=0$, $A_1(i+3) \neq \emptyset$. 
If $|A_1(i+3)|=4$, 
then there are vertices $v_j=(i+3,x_j,y_j,z_j) \in A$ $\:(j=1,2,3,4)$ 
where $x_1,x_2,x_3,x_4$ are pairwise distinct. 
Then $\{v_1,v_2,v_3,v_4\}$ is an edge of type \ref{it:102}. 
Thus, we may assume that $|A_1(i+3)| \leq 3$. 
If $A_1(i+2) = \emptyset$, then $\max_x |A_2(i+2,x)|=0$, 
so we get 
$\max_x |A_2(i+2,x)| + |A_1(i+3)| \leq 3 
= 2 + \chi[A_1(i+2) \neq \emptyset] + \chi[A_1(i+3) \neq \emptyset]$. 
Now assume that $A_1(i+2) \neq \emptyset$. 
If there is $x$ such that $|A_2(i+2,x)| \geq 4$, 
then there are vertices $v_j=(i+2,x,y_j,z_j) \in A$ $\:(j=1,2,3,4)$ 
where $y_1,y_2,y_3,y_4$ are pairwise distinct. 
Then $\{v_1,v_2,v_3,v_4\}$ is an edge of type \ref{it:102}. 
Thus, we may assume that $\max_x |A_2(i+2,x)| \leq 3$. 
We need to prove $\max_x |A_2(i+2,x)| + |A_1(i+3)| \leq 4$. 
This inequality holds if $|A_1(i+3)| \leq 1$ 
(as $\max_x |A_2(i+2,x)| \leq 3$). 
It remains to consider the case when
$\max_x |A_2(i+2,x)| = k \in\{2,3\}$ and $|A_1(i+3)| \geq 5-k$. 
We claim that in this case, 
$A$ must contain an edge of type \ref{it:103}. 
Indeed, if $k=2$, there are $x,y',y'',z',z''$ (with $y' \neq y''$) 
such that $v'=(i+2,x,y',z')$ and $v''=(i+2,x,y'',z'')$ belong to $A$. 
As $|A_1(i+3)| \geq 3$, 
there are $v_j=(i+3,x_j,y_j,z_j) \in A$ $\:(j=1,2,3)$ 
where $x_1,x_2,x_3$ are pairwise distinct. 
There are $j,k\in\{1,2,3\}$ such that $j \neq k$ and $x_j+x_k=y'+y''$. 
Then $\{v',v'',v_j,v_k\}$ is an edge of type \ref{it:103}. 
Similarly, if $k=3$, 
there are $v_j=(i+2,x,y_j,z_j) \in A$ $\:(j=1,2,3)$ 
where $y_1,y_2,y_3$ are pairwise distinct. 
As $|A_1(i+3)| \geq 2$, there are $v',v'' \in A$ where 
$v'=(i+3,x',y',z')$, $v''=(i+3,x'',y'',z'')$, and $x' \neq x''$. 
There are $j,k\in\{1,2,3\}$ such that $j \neq k$ and $y_j+y_k=x'+x''$. 
Then $\{v_j,v_k,v',v''\}$ is an edge of type \ref{it:103}. 

{\it Case 3:} $\varepsilon_i(A)=0$, $A_1(i+3) = \emptyset$. 
We need to prove 
$\max_x |A_2(i+2,x)| \leq 2 + \chi[A_1(i+2) \neq \emptyset]$. 
If $A_1(i+2) = \emptyset$ then $\max_x |A_2(i+2,x)| = 0$. 
If $A_1(i+2) \neq \emptyset$, we need to prove 
$\max_x |A_2(i+2,x)| \leq 3$. 
If there is $x$ such that $|A_2(i+2,x)| \geq 4$, 
then similarly to case {\it 2}, 
$A$ contains an edge of type \ref{it:102}. 
\end{proof}

\begin{lemma}\label{th:5dim_L2}
If $A$ is an independent set of vertices in $H_m$, then
\[
  \sum_{i\in\Z_m} \varepsilon_i(A) \; + \: 
  \sum_{i\in\Z_m} \sum_{x\in\Z_2^2} |A_2(i,x)|
  \: \leq \: 2m \: + 
  \sum_{i\in\Z_m} \chi[A_1(i) \neq \emptyset] .
\]
\end{lemma}

\begin{proof}[\bf{Proof}]
Suppose, there are $x',x''$ such that 
$|A_2(i,x')| \geq 2$, $|A_2(i,x'')| \geq 2$, 
and $x' \neq x''$. 
Then there are $v_1',v_2',v_1'',v_2'' \in A$ where 
$v_j'=(i,x',y_j',z_j')$, $v_j''=(i,x'',y_j'',z_j'')$ $\:(j=1,2)$, 
$y_1' \neq y_2'$, and $y_1'' \neq y_2''$. 
In this case, 
$\{v_1',v_2',v_1'',v_2''\}$ is an edge of type \ref{it:102}. 
Thus, for each $i\in\Z_m$, 
there is at most one $x \in \Z_2^2$ such that $|A_2(i,x)| > 1$. 
Hence, 
\[
  |A_1(i)| \: = \;
  \sum_{x\in\Z_2^2} |A_2(i,x)| 
  \: - \: \max_{x\in\Z_2^2} |A_2(i,x)|
  \: + \: \chi[A_1(i) \neq \emptyset],
\]
By \cref{th:5dim_L1}, 
\begin{multline*}
  \varepsilon_{i-3}(A) \: + \: 
  \max_x |A_2(i-1,x)| \: + \:
  |A_1(i)| \\
  \leq \: 2 \: + \: 
  \chi[A_1(i-1) \neq \emptyset] \: + \: 
  \chi[A_1(i) \neq \emptyset] ,
\end{multline*}
so we get 
\begin{multline*}
  \varepsilon_{i-3}(A) \: + \: 
  \max_x |A_2(i-1,x)| \: + \:
  \sum_{x\in\Z_2^2} |A_2(i,x)| 
  \: - \: \max_{x\in\Z_2^2} |A_2(i,x)|
  \\ \leq \: 2 \: + \: 
  \chi[A_1(i-1) \neq \emptyset] .
\end{multline*}
To finish the proof, we sum the last inequality over all $i\in\Z_m$.
\end{proof}

\begin{theorem}\label{th:5dim_T}
If $A$ is an independent set of vertices in $H_m$, then
\[
  |A| \: \leq \: 2m \: + 
  \sum_{i\in\Z_m} \chi[A_1(i) \neq \emptyset] .
\]
\end{theorem}

\begin{proof}[\bf{Proof}]
As 
$
  |A_2(i,x)| = \sum_{y\in\Z_2^2} \min\{1,\: |A_3(i,x,y)|\},
$
we get
\begin{eqnarray*}
  |A| & = & 
  \sum_{i\in\Z_m} \sum_{x,y\in\Z_2^2} |A_3(i,x,y)|
    \\ & = &
  \sum_{i\in\Z_m} \sum_{x,y\in\Z_2^2} \left(
    \max\{0,\: |A_3(i,x,y)|-1\} + \min\{1,\: |A_3(i,x,y)|\}
  \right)
    \\ & = &
  \sum_{i\in\Z_m} \varepsilon_i(A) \: + \:
  \sum_{i\in\Z_m} \sum_{x\in\Z_2^2} |A_2(i,x)| \, .
\end{eqnarray*}
Then the statement of the theorem follows from \cref{th:5dim_L2}.
\end{proof}

Let $B$ be a subset of $\Z_2^2$ of size $\lambda$, 
and $H_{m,\lambda}$ be the subgraph of $H_m$ 
induced by vertices $(i,x,y,z)$ where 
$i \in \Z_m$, $x,y \in \Z_2^2$, and $z \in B$. 
(If $\lambda=1$, then $H_{m,1}$ is isomorphic 
to the $4$-graph from Example~\ref{ex:4x4}.) 
Then ${\rm v}(H_{m,\lambda}) = 16m\lambda$.
There are 
$16 \binom{\lambda}{4}$ edges of type \ref{it:101}, 
\:$476 \lambda^4$ edges of type \ref{it:102}, 
\:$768 \lambda^4$ edges of type \ref{it:103} 
(as there exist 
$4$ choices for $x_1$, 
$\binom{4}{2}$ choices for $y_1',y_1''$, 
$2$ choices for $x_2',x_2''$, 
$4^2$ choices for $y_2',y_2''$, 
and $\lambda^4$ choices for $z_1',z_1'',z_2',z_2''$). 
There are
$120 \binom{\lambda}{2}^2$ edges of type \ref{it:104a} 
(as there exist 
$\binom{4^2}{2}$ choices for $(x_1,y_1),(x_2,y_2)$, 
$\binom{\lambda}{2}$ choices for $z_1',z_1''$, and 
$\binom{\lambda}{2}$ choices for $z_2',z_2''$). 
There are
$384 \binom{\lambda}{2} \lambda^2$ edges of type \ref{it:104b} 
(as there exist 
$4^2$ choices for $x_1,y_1$, 
$\binom{\lambda}{2}$ choices for $z_1',z_1''$, 
$4$ choices for $x_2$, 
$\binom{4}{2}$ choices for $y_2',y_2''$, and 
$\lambda^2$ choices for $z_2',z_2''$). 
There are
$512 \binom{\lambda}{2} \lambda^2$ edges of type \ref{it:104c} 
(as there exist 
$4^2$ choices for $x_1,y_1$, 
$\binom{\lambda}{2}$ choices for $z_1',z_1''$, 
$2$ choices for $x_2',x_2''$, 
$4^2$ choices for $y_2',y_2''$, and 
$\lambda^2$ choices for $z_2',z_2''$). 
Therefore, 
\begin{eqnarray*}
  {\rm e}(H_{m,\lambda}) & = &
  m \left(
    16 \binom{\lambda}{4} 
    + 1244 \lambda^4 
    + 120 \binom{\lambda}{2}^2
    + 896 \binom{\lambda}{2} \lambda^2
  \right)
    \\ & = &
  \frac{m}{3} \left(
    5168 \lambda^4 - 1536 \lambda^3 + 112 \lambda^2 - 12 \lambda
  \right) \: .
\end{eqnarray*}
When $\lambda=1$, we get ${\rm e}(H_{m,1}) = 1244m$ 
(see Example~\ref{ex:4x4}). 

For $i \in \Z_m$, 
let $V_i = \{ (i,x,y,z): \: x,y,z\in\Z_2^2 \}$. 
It follows from \cref{th:5dim_T} that 
$\alpha(H_{m,\lambda}) \leq 3m$,  
$\alpha(H_{m,\lambda} - V_i) \leq 3m-1$, 
and $\alpha(H_{m,\lambda} - (V_i \cup V_j)) \leq 3m-2$ 
if $i \neq j$. 
In particular, $V_i \cap V(H_{m,\lambda})$ 
is a critical set in $H_{m,\lambda}$ of size $16 \lambda$.
As $m \geq 4$, we can assign to each vertex $w$ of $H_{m,\lambda}$ 
two distinct critical sets which do not contain $w$. 
Now we may use \cref{th:expansion_density} with 
$c(w) = 16 \lambda$ and $d(w)=1$ for every vertex $w$ to get

\begin{corollary}\label{th:5dim_C}
If $m \geq 4$ and $\lambda \in \{1,2,3,4\}$, then
\[
  t(3m+2,4) \: \leq \:
  \left(
    24 {\rm e}
    + {\rm v} \left(
      3{\rm v}
      + 64 \lambda -\frac{43}{16}
      \right)
  \right)
  {\rm v}^{-4} ,
\]
\[
  t_*(3m+2,4) \: \leq \:
  (3m+1)^3
  \left(
    {\rm e}
    + \frac{\rm v}{24} \left(
      3{\rm v}
      + 64 \lambda -\frac{43}{16}
      \right)
  \right)
  {\rm v}^{-4} ,
\]
where ${\rm v} = 16m\lambda$ and 
$
  {\rm e} = \frac{m}{3} \left(
    5168 \lambda^4 - 1536 \lambda^3 + 112 \lambda^2 - 12 \lambda
  \right)
$.
\end{corollary}

When $m=21$, $\lambda=3$, we get 
$t_*(65,4) < 0.706335 \:$.

\section{Exact values $T(n,\alpha+1,4)$ for small ratios $n/\alpha$}
\label{sec:exact}

It is known (see \cite[Section~8]{Sidorenko:1995}) that

\begin{equation}\label{eq:401}
  T(n,\alpha+1,4) \; = \; 
  \begin{cases} 
    \;\;\;\;
    n-\alpha \;\;\;\;\;\;\;{\rm if}\;\;
      1 \leq \frac{n}{\alpha} \leq \frac{4}{3} \: , \\
    \left\lceil \frac{5}{2} n - 3 \alpha \right\rceil \;\;\;{\rm if}\;\;
	  \frac{4}{3} \leq \frac{n}{\alpha} \leq \frac{3}{2} \: .
  \end{cases}
\end{equation}

\begin{proposition}
If $\,\frac{3}{2} \leq \frac{n}{\alpha} \leq\frac{7}{4}$ 
and $n \neq \frac{7}{4}\alpha - \frac{1}{2}$, 
then 
\begin{equation}\label{eq:402}
  T(n,\alpha+1,4) \; = \; 
  \left\lceil 4n - \frac{21}{4}\alpha \right\rceil .
\end{equation}
\end{proposition}

\begin{proof}[\bf{Proof}]
The lower bound in \cref{eq:402} follows from the inequality 
${\mathrm v}(H)-\alpha(H) \leq (5{\mathrm v}(H) + 4{\mathrm e}(H)) / 21$ 
proved in \cite{Thomasse:2007}. 
To prove the upper bound, 
we will use the following inequalities: 
$T(5,4,4) \leq 5$ (trivial),
$T(6,5,4) \leq 3$ and $T(7,5,4) \leq 7$ (\cref{sec:5_16}), 
$T(8,6,4) \leq 6$ (\cref{sec:n_6_4}). 
Observe that a union of disjoint $4$-graphs $G_1,\ldots,G_m$ 
has independence number equal to the sum of the independence numbers of 
$G_1,\ldots,G_m$. 
Thus,
\[
  T(n_1+\ldots+n_m,\, \alpha_1+\ldots,\alpha_m + 1, 4) \: \leq \:
  \sum_{i=1}^m T(n_i,\alpha_i + 1,4) \, .
\]

If $\alpha = 4m$, 
then $n = 6m + k$ where $0 \leq k \leq m$. 
We get 
$
 T(6m+k,4m+1,4) \leq k T(7,5,4) + (m-k) T(6,5,4) 
 \leq 7k+3(m-k) = 4n-21m .
$

If $\alpha = 4m+1$, 
then $n = 6m + 2 + k$ where $0 \leq k \leq m-1$. 
We get 
$
  T(6m+2+k,4m+2,4) \leq T(8,6,4) + k T(7,5,4) + (m-1-k) T(6,5,4) 
  \leq 6+7k+3(m-1-k) = 4n-(21m+5) .
$

If $\alpha = 4m+2$, 
then $n = 6m + 3 + k$ where $0 \leq k \leq m-1$. 
(If $k=m$ then $n = \frac{7}{4}\alpha - \frac{1}{2}$.) 
The subcase $k=0$ is covered by \cref{eq:401}. 
When $k \geq 1$, we get 
$
  T(6m+3+k,4m+3,4) \leq 2T(8,6,4) + (k-1) T(7,5,4) + (m-1-k) T(6,5,4) 
  \leq 2 \cdot 6 + 7(k-1) + 3(m-1-k) = 4n-(21m+10) .
$

If $\alpha = 4m+3$, 
then $n = 6m + 5 + k$ where $0 \leq k \leq m$. 
We get 
$
  T(6m+5+k,4m+4,4) \leq T(5,4,4) + k T(7,5,4) + (m-k) T(6,5,4) 
  \leq 5+7k+3(m-k) = 4n-(21m+15) .
$
\end{proof}

The exceptional case $n = \frac{7}{4}\alpha - \frac{1}{2}$ 
corresponds to $n=7m+3$, $\alpha = 4m+2$. 
In this case, the right hand side of \cref{eq:402} is equal to $7m+2$, 
while we get
$T(7m+3,4m+3,4) \leq 2T(5,4,4) + (m-1) T(7,4,3) = 2 \cdot 5 + 7(m-1) = 7m+3$. 
It is likely that $T(7m+3,4m+3,4) = 7m+3$. 
This equality holds for $m=1,2$ (see \cite{Gordon:tables}).

\section{Summary of asymptotic results and open problems}\label{sec:open}

For small $k$, we tried various choices of $G_i$ 
in the circular construction 
and then used the resulting $4$-graph $H$ 
to seed the expansion construction.
We list below 
the upper bounds for the rescaled Tur\'{a}n densities $t_*(k,4)$ 
that we were able to get. 
We omit those values $k$ where the bound is weaker 
than for some $t_*(k',4)$ with $k' < k$. 
We skip proofs for $k=8,10,14$. 

\vspace{1mm}
\begin{center}
\begin{tabular}{rcl}
  $k$ & $t_*(k,4) \leq$ & reference \\
\\  4 &  1 & 
\\  5 & $5/6$ & \cref{sec:5_16}
\\  6 & $0.802611$ & Example~\ref{ex:2}
\\  7 & $0.778711$ & \cref{th:rainbow}
\\  8 & $0.765046$ & 
\\ 10 & $0.729885$ & 
\\ 14 & $0.725684$ & 
\\ 17 & $0.722438$ & Example~\ref{ex:4x4}
\\ 20 & $0.715601$ & Example~\ref{ex:4x4}
\\ 23 & $0.714739$ & Example~\ref{ex:4x4}
\\ 32 & $0.711838$ & \cref{th:5dim_C}, $m=10$, $\lambda=2$
\\ 35 & $0.709199$ & \cref{th:5dim_C}, $m=11$, $\lambda=2$
\\ 38 & $0.707575$ & \cref{th:5dim_C}, $m=12$, $\lambda=2$
\\ 41 & $0.706727$ & \cref{th:5dim_C}, $m=13$, $\lambda=2$
\\ 44 & $0.706485$ & \cref{th:5dim_C}, $m=14$, $\lambda=2$
\\ 62 & $0.706452$ & \cref{th:5dim_C}, $m=20$, $\lambda=3$
\\ 65 & $0.706335$ & \cref{th:5dim_C}, $m=21$, $\lambda=3$
\end{tabular}
\end{center}
\vspace{2mm}

Notice that our upper bound for $t_*(4)$ is not far from
$2/3$ that is the expected value of $t_*(3)$. 


\begin{problem}\label{pr4}
Is it true that $t_*(r-1) \leq t_*(r)$ ?
\end{problem}

\begin{problem}\label{pr2}
Are $t_*(r)$ bounded?
\end{problem}

\begin{problem}\label{pr3}
Are $t_*(k,r)$ nondecreasing in $k$ ?
\end{problem}

\begin{problem}\label{pr5}
Is it true that $t_*(k,r-1) \leq t_*(k,r)$ 
and $t(k,r-1) \leq \frac{k-1}{r} t(k,r)$?
\end{problem}

\section*{Acknowledgments}

The author is grateful to the two anonymous referees for their valuable suggestions.


\bibliographystyle{elsarticle-num-names-alphsort}
\bibliography{Turan}

\end{document}